# Multiple Solutions of Riemann-Type of Functional Equations

Cao-Huu T., Ghisa D., Muscutar F. A.

[1] *York University, Glendon College, Toronto, Canada*
[2] *York University, Glendon College, Toronto, Canada*
[3] *Lorain CCC, Elyria, OH* 44035, *USA*

KEY WORDS: *Riemann – type of functional equation, non trivial zeros, Riemann Hypothesis*

*Linearly independent Dirichlet L – functions satisfying the same Riemann – type of functional equation have been supposed for a long time to possess off critical line non trivial zeros. We are taking a closser look into this problem and into its connection with the Generalized Riemann Hypothesis.*

**1. Introduction**

Starting with the 1935 paper of Potter and Titchmarsh [10], the attempts to find counterexamples to the Generalized Riemann Hypothesis multiplied (see [11], [1], [2]). The idea was to use linear combinations of Dirichlet L-series which are not Euler products, but for which the analytic continuations to the whole complex plane still verify some Riemann-type of functional equations. There were two approaches to achieve this goal: one was to chose conveniently the coefficients of the respective linear combination, as in the case of the so called Davenport and Heilbronn function and the other one to use L-functions satisfying the same functional equation, for which any linear combination would do.

Potter and Titchmarsh thought they had identified two such zeros, yet they acknowledged that "the calculations are very cumbrous, and can hardly be considered conclusive". However, more zeros were indicated in [11], [1] and [2]. A mismatch of Dirichlet characters in the formula for that function brought us to think  that (see [5]) the approximation errors which naturally affected their computation produced false off critical line zeros. Our geometric function theory approach in [5] was showing that such zeros cannot exist. We discovered later our mistake and came in [4] with a correction.

In this paper we shall deal with the second approach, namely that of linearly independent L-functions satisfying the same functional equation. We perceive this class of functions not as a source of counterexamples to the GRH, but rather as a confirmation of the theory developed in [6] and [7] in which we have shown that if a general Dirichlet series can be continued analytically to the whole complex plane and the extended function satisfies a Riemann-type of functional equation then, under some other mild constraints, its nontrivial zeros have all the real part equal to 1/2.

The adopted meaning of the concept of trivial zero is that which has been originally implied, namely zeros of some elementary factors of the function (which can be trivially computed). Incidentally, for the Riemann Zeta function those zeros were the real zeros of the function, and there has rooted the opinion that the trivial zeros should be in general real. Yet for some other L-functions, as for example those defined by non primitive Dirichlet characters, it is known that trivial imaginary zeros can exist. Also, the function (1) below has zeros on the critical line, namely those of the multiplier of $\zeta(s)$ which can be trivially computed, hence they are trivial zeros. As seen in [3], only adopting this approach, the trivial zeros of the derivatives of some L-functions can be unambiguously defined.

## 2. Riemann-Type of Functional Equations Having Several Linearly Independent Solutions

The typical example of a pair of linearly independent L-functions $f_0(s)$ and $f_1(s)$ satisfying the same functional equation is that given in [1], where

$$f_0(s) = (1 + \tfrac{\sqrt{5}}{5^s})\zeta(s) \tag{2.1}$$

which coincides for $s > 1$ with the sum of the periodic Dirichlet series

$$1 + \tfrac{1}{2^s} + \tfrac{1}{3^s} + \tfrac{1}{4^s} + \tfrac{1+\sqrt{5}}{5^s} + \ldots \tag{2.2}$$

and $f_1(s)$ is the function obtained by analytic continuation to the whole complex plane of the periodic Dirichlet series

$$1 - \tfrac{1}{2^s} - \tfrac{1}{3^s} + \tfrac{1}{4^s} + \tfrac{0}{5^s} + \ldots \tag{2.3}$$

It is obvious that (2.2) and (2.3) are linearly independent, since (2.2) tends to $\infty$ as $s \to 1$, while (2.3) is convergent at $s = 1$. The functions $f_0$ and $f_1$ satisfy both the functional equation:

$$f(s) = W(s)\overline{f(1-\bar{s})}, \text{ where } W(s) = 5^{(1/2)-s}2(2\pi)^{s-1}\Gamma(1-s)\sin\tfrac{\pi s}{2} \tag{2.4}$$

Since $W(s)$ is real for real $s$, we have $\overline{f(1-\bar{s})} = f(1-s)$ and if $f(\sigma_0 + it_0) = 0$, then necessarily $f(\sigma_0 - it_0) = 0$ and due to (2.4), if $W(\sigma_0 + it_0) \neq 0$ we have also $f(1 - \sigma_0 + it_0) = 0$.

It can be easily checked that if $\overline{f_k(\bar{s})} = f_k(s)$, $k = 0, 1$ then any linear combination with real coefficients $f = \alpha_0 f_0 + \alpha_1 f_1$ of $f_0$ and $f_1$ is real on the real axis and if $f_k$ satisfies (2.4) then $f$ satisfies also (2.4). This affirmation is false if one of the coefficients $\alpha_k$ is not real.

Let us denote $\quad \varphi_\tau(s) = (1 - \tau)f_0(s) + \tau f_1(s), \quad 0 \leq \tau \leq 1 \quad$ (2.5)
where $f_0$ and $f_1$ are the functions (2.2) and (2.3)

We give to the word *deformation* used in [2] the following precise meaning: $\varphi_\tau(s)$ defined by (2.5) represents a *continuous deformation* of $f_0(s)$ into $f_1(s)$ in the sense that for any compact set $K \subset \mathbb{C} \setminus \{1\}$ we have $\lim_{\tau \to \tau'} \varphi_\tau(s) = \varphi_{\tau'}(s)$ uniformly in $K$ and $\varphi_0(s) = f_0(s)$ respectively $\varphi_1(s) = f_1(s)$.

**Theorem 1**: *For every $\tau \in [0, 1]$, the function $\varphi_\tau(s)$ given by (2.5) satisfies (2.4) and as $\tau$ varies from $0$ to $1$, $\varphi_\tau(s)$ represents a continuous deformation of $f_0(s)$ into $f_1(s)$*

*Proof:* Since $f_0(s)$ and $f_1(s)$ satisfy (2.4), and $f_k(\bar{s}) = \overline{f_k(s)}$ we have:
$\varphi_\tau(s) = (1 - \tau)f_0(s) + \tau f_1(s) =$
$(1 - \tau)W(s)f_0(1 - s) + \tau W(s)f_1(1 - s) = W(s)[(1 - \tau)f_0(1 - s) + \tau f_1(1 - s)] = W(s)\varphi_\tau(1 - s)$,
hence $\varphi_\tau(s)$ satisfies (2.4). When $\tau = 0$ we get $\varphi_0(s) = f_0(s)$ and when $\tau = 1$ we get $\varphi_1(s) = f_1(s)$.

For any compact $K$, $\max_k \sup_{s \in K} |f_k(s)|$ is a finite number $M$, due to the uniform continuity of $f_k$ on $K$. Then $|\varphi_\tau(s) - \varphi_{\tau'}(s)| = |\tau' - \tau||f_0(s) - f_1(s)| \leq 2M|\tau' - \tau|$, which proves the theorem.¤

Although there were no specific values indicated, it has been implied in [1] that some of the zeros of $\varphi_\tau(s)$ should be off critical line. However, we can prove the following:

**Theorem 2**. *For every $\tau \in [0, 1]$ the non trivial zeros of $\varphi_\tau$ have the real part equal to $1/2$.*

*Proof:* Let us notice that $f_0(s)$ has the same non trivial zeros as $\zeta(s)$ and $f_1(s)$ is the Dirichlet L-function $L(5, 3, s)$. By the Generalized Riemann Hypothesis, which we are taken as true in this paper, the non trivial zeros of the two functions have the real part $1/2$. To prove the theorem we need only to check that $\varphi_\tau$ is the analytic continuation to the whole complex plane, except for the pole $s = 1$ of the corresponding linear combination of the two series, which is obvious, and that $\lim_{\sigma \to +\infty} \varphi_\tau(\sigma + it) = 1$. Indeed, this last equality is true having in view that both $f_0(s)$ and $f_1(s)$ have the limit $1$ as $\sigma \to +\infty$. Since (2.1) and (2.2) have both the abscissa of convergence $\sigma_c = 0$, we have also that for every $\tau$, $0 \leq \tau \leq 1$, the abscissa of convergence of the respective linear combination of the two series is $0$. Then, by [8], Theorem 3, all the non trivial zeros of $\varphi_\tau(s)$ have the real part $1/2$.¤

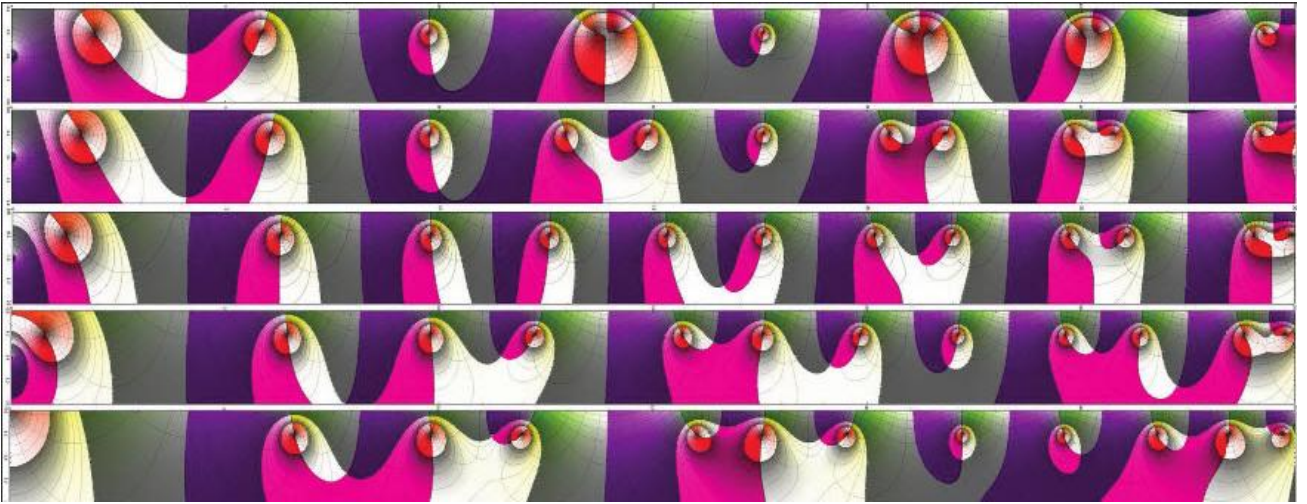

Fig. 1.   $\operatorname{Re}\varphi_\tau(s) = 1/2$ for every $\tau$,  $0 \leq \tau \leq 1$

Fig. 1 illustrates this situation for the particular values of $\tau$ : $0$, $0.25$, $0.5$, $0.75$, $1$ and $s \in [0,1]x[0,30]$. It can be seen how the zeros of $f_0(s)$ evolve alongside the critical line into those of $f_1(s)$ as $\tau$ varies from $0$ to $1$. For space economy, the axes are rotated.

**Theorem 3**. *Suppose that two L-functions $f_0$ and $f_1$ are continuations of Dirichlet series having the abscissa of convergence less than $1/2$ and they are such that* $\lim_{\sigma \to +\infty} f_k(\sigma + it) = 1$, $k = 0,1$. *If the two functions satisfy the same Riemann-type of functional equation, then any interval $I = \{1/2 + it \mid t_1 \leq t \leq t_2\}$ of the critical line can be extended to an interval $I'$ such that $f_0$ and $f_1$ have the same number of zeros in $I'$.*

*Proof:* We notice that $\varphi_\tau(s)$ defined by $(2.5)$ satisfies also the respective functional equation and the Dirichlet series of $\varphi_\tau(s)$ has the abscissa of convergence less than $1/2$. Then, by [8], Theorem 3, $\varphi_\tau(s)$ satisfies the RH. Let us follow the trajectory of one particular zero of $f_0(s)$ when $\tau$ varies from $0$ to $1$. We have $f_0(s) - \varphi_\tau(s) = \tau[f_0(s) - f_1(s)]$. In particular, if $f_0(s_0) = 0$, then $|\varphi_\tau(s_0)| = \tau|f_1(s_0)|$ and for every $\epsilon > 0$ there is $\delta > 0$ such that $|\varphi_\tau(s_0)| < \epsilon$ if $\tau < \delta$, which means that the disc centered at $\varphi_\tau(s_0)$ and of radius $2\epsilon$ contains the origin. Thus the pre-image by $\varphi_\tau(s)$ of this disc contains a zero of $\varphi_\tau(s)$. Let $s_\tau$ be the closest zero to $s_0$ as previously defined. We can interpret $\varphi_\tau(s)$ as a deformation of $f_0(s)$ which *carried* with it the zero $s_0$ into a new location $s_\tau$ on the critical line. At the next step, the deformation of $\varphi_\tau(s)$, which can be checked to be a new deformation of $f_0(s)$ will carry $s_\tau$ to a new location $s_{\tau'}$ and so on. At every step $s_0$ is moved in the same direction into an open neighborhood on $I$ and since $I$ is a compact set, after a finite number of steps $s_0$ is carried to a zero $s_1$ of $f_1(s)$. This process establishes a one to one correspondence between the zeros of $f_0(s)$ and those of $f_1(s)$ such that if an interval $I$ for $t$ is given we can extend it if necessary to a bigger one which contains only pairs of corresponding zeros of the two functions and the theorem is completely proved ¤

*Remarks*: 1. We notice that the trajectories of the zeros defined in this theorem cannot intersect the boundaries of the strips $S_k$, respectively $S'_k$ corresponding to the two functions, since there are no zeros on those boundaries. This means that we can take $I'$ as being the intersection of the critical line with $S_k \cap S'_k$. In other words, in every strip $S_k$ defined by one of the two L-functions satisfying the same functional equation, both of them have the same number of zeros. In particular, the functions $N(T)$ counting the non trivial zeros of the two functions in the interval $[0, T]$ have values different by at most one unit.

2. An interesting question arises about some of the non trivial zeros of the function $f_1(s)$ obtained by analytic continuation of the Dirichlet series $(2.3)$. It can be noticed that the zeros of $1 + 5^{(1/2)-s}$ are preserved during the continuous deformation of $f_0(s)$ into $f_1(s)$, in other words some of the non trivial zeros of $f_1(s)$ are the same as some trivial zeros of $f_0(s)$. This apparent contradiction can be settled by noticing that $f_1(s)$ can be also factorized by $1 + 5^{(1/2)-s}$. Indeed $1 - \frac{1}{2^s} - \frac{1}{3^s} + \frac{1}{4^s} + \frac{0}{5^s} + \frac{1}{6^s} - \frac{1}{7^s} - \frac{1}{8^s} + \frac{1}{9^s} + \frac{0}{10^s} + \ldots =$
$(1 + \frac{\sqrt{5}}{5^s})(1 - \frac{1}{2^s} - \frac{1}{3^s} + \frac{1}{4^s} - \frac{\sqrt{5}}{5^s} + \frac{1}{6^s} - \frac{1}{7^s} - \frac{1}{8^s} + \frac{\sqrt{5}}{10^s} + \ldots)$. This last series has the same abscissa of convergence as the series $(2.3)$ and it can be continued analytically to the whole complex plane to an L-function $f_2(s)$. Hence the true non trivial zeros of $f_1(s)$ are only those which are non trivial also for $f_2(s)$.

## 3. More Dirichlet L-Functions Verifying the Same Functional Equation

Let $L(s; \chi_k)$ be the Dirichlet L-function defined by a Dirichlet character $\chi_k$ modulo $q$, i.e.

$$L(s; \chi_k) = \Sigma_{n=1}^{\infty} \chi_k(n) / n^s, \qquad (3.1)$$

where $\chi_k(n)$ is the $k$-th Dirichlet character modulo $q$. Then $L(s; \chi_k)$ satisfies (see [9]) the functional equation

$$L(s; \chi) = \epsilon(\chi) W(s) L(1 - s; \overline{\chi}), \qquad (3.2)$$

where $W(s) = 2^s q^{(1/2)-s} \pi^{s-1} \Gamma(1-s) \sin \frac{\pi}{2}(s + \kappa)$ and $\kappa = 0$ if $\chi(-1) = 1$ ($\chi$ is even), $\kappa = 1$ if $\chi(-1) = -1$ ($\chi$ is odd), $\epsilon(\chi) = \tau(\chi)/i^\kappa \sqrt{q}$ and $\tau(\chi) = \Sigma_{k=1}^{q} \chi(k) Exp\{2\pi i/q\}$, $|\epsilon(\chi)| = 1$

Obviously, if $\chi_k$ and $\chi_{k'}$ have different parities, the corresponding functional equations are different and if $\chi_k$ and $\chi_{k'}$ have the same parity in order for them to coincide, we must have $\tau(\chi_k) = \tau(\chi_{k'})$. In other words, $L(s; \chi_k)$ a nd $L(s; \chi_{k'})$ satisfy the same functional equation if and only if $\tau(\chi_k) = \tau(\chi_{k'})$. Yet, the line matrices representing the values of different $\chi_k$ are linearly independent and therefore $k \neq k'$ implies $\tau(\chi_k) \neq \tau(\chi_{k'})$. Therefore, no two Dirichlet L-functions can satisfy the same Riemann-type of functional equation. In order to find L-functions satisfying the same functional equation we need to expand on the example from th e section 2, considering an arbitrary modulus $q$. So, let

$$f_0(s) = [1 + q^{(1/2)-s}]\zeta(s) \qquad (3.3)$$

and let us notice that $f_0(s)$ satisfies the functional equation

$$f(s) = W(s)f(1-s), \text{ where } W(s) = q^{s-1/2}2(2\pi)^{s-1}\Gamma(1-s)\sin\tfrac{\pi s}{2} \qquad (3.4)$$

Then we should look for even real primitive non principal characters $\chi(\mod q)$, which generate Dirichlet L-functions satisfying the functional equation (3.4). The first such candidate is $\chi_3(\mod 5)$, which we have already encountered. Then, ordered by increasing $q$ the next characters of interest are: $\chi_2(\mod 8)$, $\chi_3(\mod 10)$, $\chi_4(\mod 12)$, $\chi_7(\mod 13)$, $\chi_3(\mod 15)$, $\chi_9(\mod 17)$, $\chi_{10}(\mod 21)$, etc. We cannot find a formula generating all these characters, but is reasonable to assume that their sequence is infinite.

There are also odd quadratic characters satisfying a similar equation, namely $\chi_2(\mod 3)$, $\chi_2(\mod 4)$, $\chi_2(\mod 6)$, $\chi_4(\mod 7)$, $\chi_6(\mod 11)$, $\chi_3(\mod 12)$, $\chi_4(\mod 14)$, $\chi_5(\mod 15)$, $\chi_{10}(\mod 19)$, etc. We checked that for all these characters, even and odd, $\epsilon(\chi) = 1$, but since $\kappa = 1$ for odd characters, we have $\cos\tfrac{\pi s}{2}$ instead of $\sin\tfrac{\pi s}{2}$ in the expression of $W(s)$.

The Dirichlet L-functions defined by even characters satisfy each one the equation (3.4) corresponding to the respective $q$. To find a similar property for the odd caracters, we need to replace $\zeta(s)$ in the formula (3.3) with a Davenport and Heilbronn function modulo $q$ (see [4]), which contains in (3.4) the right trigonometric function. Hence instead of (3.3) we should take:

$$f_0(s) = \tfrac{1}{2}\{[L(s;\chi) + L(s;\overline{\chi})] + i\tan\theta[L(s;\chi) - L(s;\overline{\chi})]\} \qquad (3.5)$$

where $\chi$ and $\overline{\chi}$ are complex conjugate Dirichlet characters modulo $q$ and $\epsilon(\chi) = e^{2i\theta}$

The Dirichlet L-functions defined by these Dirichlet characters satisfy each one the equation (3.4) corresponding to the respective $q$, since in (3.2), $\kappa = 0$ and it can be easily checked that all have $\epsilon(\chi) = 1$.

**Theorem 4**: *For any couple $q$ and $k$ from the previous sequence the function $L(q,k,s)$ has some trivial zeros on the critical line, namely the zeros of $1 + q^{(1/2)-s}$ and its non trivial zeros are obtained by moving to their location the zeros of $f_0(s)$ from (3.3), respectively (3.5) following trajectories on the critical line by a continuous deformation of $f_0(s)$ into $L(q,k,s)$. In particular, every interval $I$ of the critical line can be extended to an interval $I'$ such that the number of non trivial zeros situated in $I'$ of $L(q,k,s)$ and of $f_0(s)$ is the same.*

The proof is similar to that of Theorem 3 and we shall skip it.